\documentclass[12pt,reqno]{amsart}
\usepackage{hyperref}
\usepackage{amsmath,amssymb, amsthm}
\usepackage{caption}
\usepackage{amssymb,amsmath,euscript,enumerate,tikz}
\usepackage[margin=1in]{geometry}
\usepackage{graphicx}
\usepackage{float}
\usepackage{pgf,tikz}
\usepackage{multirow}
\usepackage{mathrsfs}
\usepackage{subfigure}
\usepackage{pgfplots}
\pgfplotsset{compat=1.15}
\usepackage{mathrsfs}
\usetikzlibrary{arrows}
\usepackage{mathtools}

\DeclareMathOperator{\reg}{reg}

\newcommand \iv{\operatorname{iv}}

\newcommand \cdeg{\operatorname{cdeg}}

\newcommand \Tor{\operatorname{Tor}}

\newcommand \dist{\operatorname{dist}}

\newcommand{\J}{\mathcal{J}}

\theoremstyle{plain}
\newtheorem{theorem}{Theorem}[section]
\newtheorem{lemma}[theorem]{Lemma}

\newtheorem{corollary}[theorem]{Corollary}

\theoremstyle{definition}

\newtheorem{notation}[theorem]{Notation}

\newtheorem{remark}[theorem]{Remark}
\newtheorem{definition}[theorem]{Definition}
\newtheorem{example}[theorem]{Example}

\begin{document}
\title[Application of Betti Splittings to the Regularity of Binomial Edge Ideals]{Application of Betti Splittings to the Regularity of Binomial Edge Ideals}
\author[Rajiv Kumar]{Rajiv Kumar}
\email{rajiv.kumar@iitjammu.ac.in}
\author[Paramhans Kushwaha]{Paramhans Kushwaha}
\email{2022rma2004@iitjammu.ac.in}
\address{Indian Institute of Technology, Jammu}
\address{NH-44, PO Nagrota, Jagti, Jammu and Kashmir 181221}
\date{\today}


	\begin{abstract}
        Let $G$ be a simple graph on $[n]$ vertices and $J_G$ denote the corresponding binomial edge ideal in $S=k[x_1,\ldots,x_n,y_1,\ldots,y_n]$, where $k$ is a field. In this paper, we provide an upper bound for the regularity of binomial edge ideals of trees. Our approach leverages the recently developed framework of Betti splittings of binomial edge ideals, a powerful technique for studying homological invariants via decompositions into simpler ideals. We also give a lower bound for the same. As a consequence, we explicitly compute the regularity of binomial edge ideals of trees under certain conditions.
	\end{abstract}
	
	\maketitle
  \section{introduction}

Let $G$ be a finite simple graph with vertex set $[n]:=\{1,2, \ldots,n\}$ and edge set $E(G)$. Consider $S=k[x_1,\ldots,x_n,y_1,\ldots, y_n]$ a polynomial ring over a field $k$. Herzog et al. \cite{BEIandConditionalDependance} and independently Ohtani \cite{GraphsandIdealsGeneratedbysome2-minor-Ohtani-Masahiro} introduced a notion of binomial edge ideal corresponding to a finite simple graph. The \emph{binomial edge ideal} of $G$, denoted $J_G$, is defined as $$J_G:=\langle f_{ij}~|~\{i,j\}\in E(G)\rangle\subset S,\text{ where $f_{ij}=x_iy_j-x_jy_i$.}$$
Combinatorially, $J_G$ can be interpreted as the ideal generated by a subset of maximal minors of a $2\times n$ generic matrix, corresponding precisely to the edges of $G$. This framework simultaneously generalizes two classical constructions: the determinantal ideal and the ideal of adjacent minors of a $2\times n$ generic matrix. Various authors have explored the algebraic and homological properties of binomial edge ideals in terms of the combinatorial properties of the underlying graph $G$, e.g., (see \cite{CM-BEI-Ene-H-H},\cite{HKS17},\cite{OntheextremalBettinumofBEIofblockgraphds},\cite{HilbertFunctionof-BEI},\cite{CBEIPeeva}) and their references. Notably, Herzog et al. \cite{BEIandConditionalDependance} described the reduced Gr\"obner basis of $J_G$ in terms of admissible paths of $G$, and proved that all binomial edge ideals $J_G$ are radical. Further, the minimal primes of $J_G$ are characterised combinatorially in terms of minimal cut sets of $G$. 

The study of Betti numbers and Castelnuovo–Mumford regularity of homogeneous ideals has gained significant interest in recent years, due to their deep connections with algebraic and geometric properties. For a finite graded $S$-module $M$, the \emph{Castelnuovo-Mumford regularity} (henceforth called \emph{regularity}) of $M$ is $\reg_S(M):=\sup\{j-i:\beta_{i,j}(M)\neq 0\}$, where $\beta_{i,j}(M):=\dim_k(\Tor_i^S(k,M)_j)$ is the $(i,j)$-th \emph{graded Betti number} of $M$.  For brevity, we write $\reg(M)$ instead of $\reg_S(M)$ when the ring is clear from the context.

 The main aim of this paper is to understand the regularity of binomial edge ideals of trees. Our approach employs the recently developed technique of Betti splitting introduced in \cite{Splittingsofmonomialideals}, a powerful tool to study homological properties. Matsuda and Murai \cite{RegboundsforBEI-MM} proved that if $G$ is a graph on the vertex set $[n]$, then $l\leq \reg(S/J_G)\leq n-1$, where $l$ is the length of a longest induced path in $G$. Further, they conjectured that $\reg(S/J_G)=n-1$ if and only if $G$ is a path of length $n-1$. This conjecture was settled in \cite{TheCM-regof-BEI-Kiani-Madani} by Kiani and Madani. They also conjectured in \cite{On-the-BEI-of-pair-of-graphs-Madani-Kiani} that
$\reg(S/J_G)\leq c(G),$
where $c(G)$ denotes the number of maximal cliques of $G$. This conjecture was subsequently proved by the authors of \cite{A-proof-for-a-conjecture-reg-BEI-Malayeri-Madani-Kiani}, who established a new combinatorial upper bound (see \cite[Theorem 2.4 and Corollary 2.7]{A-proof-for-a-conjecture-reg-BEI-Malayeri-Madani-Kiani}). For a broader overview of developments on this topic, we refer the reader to the recent survey \cite{jayanthan-Kumar-2025generalizedbinomialedgeideals} and the references therein.

  Jayanthan, Narayanan, and Rao \cite{RegBEIcertainblockgraph} proved that for a tree $T$, $\reg(S/J_T)\geq \iv(T)+1$ with equality if and only if $T$ is jewel-free, where $\iv(T)$ is the number of internal vertices of $T$.  A \emph{jewel} $\mathcal{J}$ is a tree obtained from three copies of $K_{1,3}$ by identifying one free vertex from each. The existence of a jewel in a tree significantly influences the regularity of the corresponding binomial edge ideal. For a tree $T$, authors give an upper bound for $\reg(S/J_T)$ \cite[Corollary 0.1]{Anupperbound-regboundof-trees-Jayantha-Narayanan-Rao}. Another notable upper bound, which is currently the best known for binomial edge ideals of trees, is given by Kumar in \cite[Theorem 4.5]{BEI-of--generalized-block-graphs-Arvind-Kumar}. Mascia and Rinaldo \cite[Corollary 3.2]{Krull-dim-and-reg-BEI-of-Block-graphs-Rinaldo} have established the strongest known lower bound for the regularity of binomial edge ideals of trees.

Finding the exact regularity of binomial edge ideals in terms of the properties of the underlying graph is a challenging problem.
Betti splitting is one of the important tools for studying the regularity of certain ideals. Jayanthan, Sivakumar, and Van Tuyl recently described the Betti splittings of binomial edge ideals \cite[Theorem 4.2, Theorem 4.3]{Partial-Betti-splittings-applications-to-BEI-Jayanthan-Shivakumar-Van-tuyl}. In this article, our goal is to determine the regularity of binomial edge ideals for trees using Betti splitting. We provide upper and lower bounds for the same as follows: \vspace {2mm}

\noindent
\textbf{Theorem~~\ref{upper bound for reg}.}
    Let $G$ be a tree containing $s$ generalized jewels. Suppose $G_1,\ldots,G_p$ are connected components of $G_{\mathcal{J}}$, and $|V(G_i)|=s_i$ for all $1\leq i\leq p$. Denote $e_G$ the number of isolated edges in $G_{\mathcal{J}}$. Then
        $$\reg(S/J_G)\leq \iv(G)+1+D_G-2s-\sum\limits_{i=1}^{p}\left\lfloor \frac{s_i}{3}\right\rfloor-e_G.$$
\vspace{2mm}

\noindent
\textbf{Theorem~~\ref{lower bound for reg}.}
   Under the hypothesis of Theorem \ref{upper bound for reg}, we have  $$\iv(G)+1+D_G-3s-\mu_G+p\leq\reg(S/J_G).$$
 Corollary \ref{Trees of exact reg} gives the exact regularity of trees under certain conditions, which shows both the bounds in Theorem \ref{upper bound for reg} and Theorem \ref{lower bound for reg} are sharp. Example \ref{Exm bound better than kumar} shows that the gap between the upper bound in Theorem \ref{upper bound for reg} and the one from \cite{BEI-of--generalized-block-graphs-Arvind-Kumar} for the regularity of the binomial edge ideal for trees can be made arbitrarily large. The lower bound in \cite{Krull-dim-and-reg-BEI-of-Block-graphs-Rinaldo} for the regularity of binomial edge ideals of trees is always smaller than that of Theorem \ref{lower bound for reg} if there are two centers $v, w$ of jewels with $D(v), D(w)\geq 4$. 

Our article is organized in the following manner. Section \ref{SecPre} contains all the necessary definitions and literature that are required for the article. In Section \ref{SecBettiSpilliting}, we use the Betti splitting to develop the results that are helpful to study the regularity of binomial edge ideals of trees. Section \ref{SecMainResults} is devoted to proving the main results of this article, see Theorem \ref{upper bound for reg}, \ref{lower bound for reg}, and Corollary \ref{Trees of exact reg}.
\subsection*{Acknowledgments}
The authors would like to thank referees for their careful reading and many helpful suggestions that have improved the readability and overall quality of the paper. The first author is supported by the Indian Institute of Technology Jammu, under the Seed Grant (SGR-100014).
\section{preliminaries}\label{SecPre}
    In this section, we review the essential background in graph theory and commutative algebra that is necessary for the later results.
    \subsection{Background from Graph Theory} 
        For a graph $G$, denote $V(G)$ the vertex set, and $E(G)$ the edge set of $G$. We say a graph $G$ is \emph{simple} if it does not have multiple edges or loops. A graph $H$ is called a \emph{subgraph} of $G$ if $V(H)\subset V(G)$ and $E(H)\subset E(G)$. A subgraph $H$ is an \emph{induced subgraph} of $G$ on $V(H)$, if for all $i,j\in V(H)$, $\{i,j\}\in E(G)$ if and only if $\{i,j\}\in E(H)$. For $v\in V(G)$, the \emph{open neighborhood} of $v$ is $N_G(v)=\{u\in V(G):\{u,v\}\in E(G)\}$, and the \emph{closed neighborhood} of $v$ is $N_G[v]=N_G(v)\cup \{v\}$. A vertex $u \in V(G)$ is \emph{adjacent} to a vertex $v$ if $u\in N(v).$ The \emph{degree} of a vertex $v\in V(G)$ is $\deg_G(v)=|N_G(v)|$. A vertex $v$ is said to be a \emph{leaf} or \emph{pendant} vertex if $\deg_G(v)=1$. A \emph{complete graph} or a \emph{clique} $K_n$ is a graph on $[n]$ such that $\{i,j\}\in E(G)$ for all $i\neq j\in [n]$. An induced subgraph $C$ is a \emph{clique of $G$} if it is a complete graph. A \emph{maximal clique} of $G$ is a clique that is not an induced subgraph of any other clique of $G$. For $v\in V(G)$, let $\cdeg_G(v)$ denote the number of the maximal cliques that contain $v$, called the \emph{clique degree} of $v$. A vertex of clique degree one is called a \emph{free} vertex, e.g., a leaf vertex is a free vertex. On the other hand, a vertex $v\in V(G)$ is called an \emph{internal vertex} if $\cdeg_G(v)\geq 2$. The number of internal vertices of $G$ is denoted by $\iv(G)$. A \emph{path} $P$ from $v$ to $w$ of length $l$ in $G$ is a sequence $P:v=v_0,\ldots,v_l=w$ of non-repeating vertices such that $\{v_i,v_{i+1}\}\in E(G)$ for all $0\leq i\leq l-1$. The \emph{distance} between $v$ and $w$ is the length of the shortest path from $v$ to $w$, denoted $\dist(v,w)$. A graph $G$ is a \emph{tree} if for any $v\neq w\in V(G)$, there is a unique path from $v$ to $w$. A longest path in a tree is called a \emph{spine} of the tree.
\begin{definition}
    Let G be a simple graph on the vertex set $[n]$, $e=\{i,j\}\notin E(G)$, and $e'=\{i',j'\}\in E(G)$. The simple graphs $G\setminus {e'}$, $G\cup \{e\}$, and $G_e$ on the vertex set $[n]$ are defined with the following edge sets:
    \begin{enumerate}
        \item $E(G\setminus {e'})=E(G)\setminus \{e'\}$;
        \item $E(G\cup\{e\})=E(G)\cup \{e\}$;
        \item $E(G_e)=E(G)\cup \{\{k,l\}:k,l\in N(i)\text{ or }k,l\in N(j)\}$.
        \item For $v\in [n]$, the simple graph $G\setminus{\{v\}}$ or $G\setminus{v}$ on the vertex set $[n]\setminus{\{v\}}$ is defined with edge set $E(G\setminus{v})=E(G)\setminus{\{\{u,v\}: u\in N(v)\}}$.    
    \end{enumerate}
\end{definition}
A graph $G$ is a \emph{clique-sum} of $G_1$ and $G_2$ if $V(G_1)\cup V(G_2)=V(G)$, $E(G_1)\cup E(G_2)=E(G)$, and the induced subgraph of $G$ on $V(G_1)\cap V(G_2)=K_t$, and we write $G=G_1\cup_{K_t} G_2$. If $t=1$, then we write $G_1\cup _{v} G_2$ for the clique-sum $G_1\cup_{K_1} G_2$, where $V(K_1)=\{v\}$. We say a graph $G$ is \emph{decomposable} if there exist induced subgraphs $G_1$ and $G_2$ of $G$ such that $G=G_1\cup _{v} G_2$ with $v$ a free vertex in both $G_1$ and $G_2$, and we say $G=G_1\cup G_2$ is a decomposition of $G$. 

The following result gives the behavior of regularity of binomial edge ideals under the decomposition of a graph and is proved in [\cite[Proposition 3]{OntheextremalBettinumofBEIofblockgraphds}, \cite[Theorem 3.1]{RegBEIcertainblockgraph}].
\begin{theorem}\label{Gluing theorem}
    Let $G=G_1\cup G_2$ be a decomposition of $G$. Then $$\reg(S/J_G)=\reg(S/J_{G_1})+\reg(S/J_{G_2}).$$
\end{theorem}

We recall the following result, due to Matsuda and Murai, which is frequently used throughout the article.
\begin{corollary}\label{reg reln b/w inducedsubgraph}\cite[Corollary 2.2]{RegboundsforBEI-MM}
    Let $H$ be an induced subgraph of $G$. Then $$\reg(S/J_G)\geq \reg(S/J_H).$$
\end{corollary}
We describe the construction of a useful short exact sequence introduced by Ene, Herzog, and Hibi \cite{CM-BEI-Ene-H-H}.
\subsection{Ene-Herzog-Hibi (EHH) process }
Let $G$ be a block graph, $\Delta(G)$ be the clique complex of $G$, and $F_1,\ldots, F_r$ be a leaf order on the facets of $\Delta(G)$. Assume that $r>1$. Let $v\in G$ be the unique vertex in $F_r$ such that $F_r\cap F_j\subseteq\{v\}$ for all $j<r$. Let $G'$ be the graph obtained by adding the necessary edges to $G$ so that $N[v]$ is a clique. Let $G''$ be the graph induced on $G\setminus\{{v}\}$ and $\Tilde{G}$ be the graph induced on $G'\setminus{\{v\}}$. Then there exists an exact sequence \begin{equation}\label{EHH-process}
    0\rightarrow{S/J_G}\rightarrow S/J_{G'}\oplus S/J_{G''}\rightarrow S/J_{\Tilde{G}}\rightarrow 0.
\end{equation}
We call $G', G''$, and $\Tilde{G}$ the graphs obtained by applying the EHH-process on $G$ with respect to $v$. 

\subsection{Betti splittings of binomial edge ideals}
 The notion of Betti splitting was first introduced by Francisco, Tai Ha, and Van Tuyl \cite{Splittingsofmonomialideals} for the monomial ideals. As an application, they showed that edge ideals of graphs are splittable, and computed the Betti numbers of the cover ideals of Cohen-Macaulay bipartite graphs. Recently, Jayanthan, Sivakumar, and Van Tuyl \cite{Partial-Betti-splittings-applications-to-BEI-Jayanthan-Shivakumar-Van-tuyl} studied the Betti splittings of homogeneous ideals. Let $S=k[x_1,\ldots,x_n]$ be a polynomial ring over a field $k$. For a homogeneous ideal $L$ in $S$, denote $\mathfrak{G}(L)$ a set of minimal generators of $L$. Given a homogeneous ideal $I$, split this ideal as $I = J + K$, where $\mathfrak{G}(I)$ is the disjoint union of $\mathfrak{G}(J)$ and $\mathfrak{G}(K)$. 

\begin{definition}\cite[Definition 3.4]{Partial-Betti-splittings-applications-to-BEI-Jayanthan-Shivakumar-Van-tuyl}
   Let $I,J$ and $K$ be homogeneous ideals of $S$ with respect to the standard $\mathbb N$-grading such that $\mathfrak{G}(I)$ is the disjoint union of $\mathfrak{G}(J)$ and $\mathfrak{G}(K)$. Then \( I = J + K \) is a \emph{Betti splitting} if
\[
\beta_{i,j}(I) = \beta_{i,j}(J) + \beta_{i,j}(K) + \beta_{i-1,j}(J \cap K)
\]
for all \( i ,j\geq 0\).
\end{definition}
Once $I=J+K$ is a Betti splitting, then one can extract information such as regularity, projective dimension, and linearity defect of $I$ from those of $ J, K$, and $J\cap K$. Authors in \cite{Partial-Betti-splittings-applications-to-BEI-Jayanthan-Shivakumar-Van-tuyl} established Betti splitting of binomial edge ideals, and proved the following result:
\begin{corollary}\label{reg BEI for betti splitting}\cite[Corollary 4.4]{Partial-Betti-splittings-applications-to-BEI-Jayanthan-Shivakumar-Van-tuyl} Let $e = \{{u, v}\} \in E(G)$ with $v$ a leaf vertex and suppose $e$ is not the only edge of $G$. Then $$\reg(S/J_G)=\max\{\reg(S/J_{G\setminus{e}}),\reg(S/J_{(G\setminus{e})_e})+1\}.$$
\end{corollary}

\section{application of betti splitting to the regularity}\label{SecBettiSpilliting}
 In this section, we establish several technical results to provide an upper bound on the regularity of binomial edge ideals for trees.
 
The graph on the vertex set $[n]\cup\{v\}$ with edge set $\{\{v,i\}:1\leq i\leq n\}$ is called a \emph{star} graph, denoted $K_{1,n}$, and the vertex $v$ is called the \emph{center} vertex of the star graph. For a graph $G$, $v\in V(G)$, and $2\leq s\in \mathbb N$, denote $G_s^v$ and $\overline{G}_s^v$ graphs obtained by identifying $v$ with a vertex of $K_s$, i.e., $G_s^v=G\cup _vK_s$ and a leaf vertex of $K_{1,s}$, i.e., $\overline{G}_s^v=G\cup _vK_{1,s}$, respectively. The following figure illustrates the construction of $G_s^v$ and $\overline{G}_{s,t}^v$.
\begin{center}
\begin{tikzpicture}[scale=.35, line cap=round,line join=round,x=1.0cm,y=1.0cm]
\clip(-2.585681940000006,-4.946980480000011) rectangle (20.4921661000000075,4.358366320000012);
\draw  (0,0)-- (-2,0);
\draw  (0,0)-- (0,2);
\draw (0,0)-- (0,-2);
\draw  (6,0)-- (4,0);
\draw  (6,0)-- (6,2);The following example illustrates Lemmas 3.1 and 3.2 and shows how they can be effectively applied to study the regularity of trees.
\draw  (6,0)-- (6,-2);
\draw  (6,0)-- (8,0);
\draw  (8,0)-- (8,2);
\draw  (8,0)-- (8,-2);
\draw  (8,0)-- (10,0);
\draw  (16,0)-- (14,0);
\draw  (16,0)-- (16,2);
\draw  (16,0)-- (16,-2);
\draw  (16,0)-- (18,-2);
\draw  (16,0)-- (18,2);
\draw  (18,2)-- (20,0);
\draw  (18,-2)-- (20,0);
\draw  (16,0)-- (20,0);
\draw  (18,2)-- (18,-2);
 \begin{scriptsize}
 \fill  (0,0) circle (3.5pt);
 \fill (-2,0) circle (3.5pt);
 \fill  (0,2) circle (3.5pt);
 \fill  (0,-2) circle (3.5pt);
 \fill  (6,0) circle (3.5pt);
 \fill  (4,0) circle (3.5pt);
 \fill  (6,2) circle (3.5pt);
 \fill  (6,-2) circle (3.5pt);
 \fill  (8,0) circle (3.5pt);
 \fill  (8,2) circle (3.5pt);
 \fill  (8,-2) circle (3.5pt);
 \fill  (10,0) circle (3.5pt);
 \fill (16,0) circle (3.5pt);
 \fill  (14,0) circle (3.5pt);
 \fill  (16,2) circle (3.5pt);
 \fill  (16,-2) circle (3.5pt);
 \fill  (18,-2) circle (3.5pt);
 \fill  (18,2) circle (3.5pt);
 \fill  (20,0) circle (3.5pt);
\end{scriptsize}
\draw[color=black] (-.4,-.35) node {$v$};
\draw[color=black] (15.6,-.35) node {$v$};
\draw[color=black] (5.6,-.35) node {$v$};
\draw[color=black] (0,-3.5) node {$G$};
\draw[color=black] (7,-3.5) node {$\overline{G}_4^v$};
\draw[color=black] (17,-3.5) node {$G_4^v$};
\end{tikzpicture}
\end{center}
\begin{lemma}\label{gluing Ks and star}
    Adopt the notation as above. Then $\reg(S/J_{\overline{G}_s^v})=1+\reg(S/J_{G_s^v})$ for all $s\geq 2$.
\end{lemma}
\begin{proof}
   We apply induction on $s$ to prove the lemma. Let $v_1,\ldots,v_s$ be leaves and $c$ be the center vertex of $K_{1,s}$. Suppose $v=v_1$ in $\overline{G}_s^v$. For $s=2$, notice that $\overline{G}_2^v=\{v_2,c\}\cup G_2^v$ is a decomposition of $G$. Therefore, the assertion follows from Theorem \ref{Gluing theorem}. Assume $s\geq 3$. Observe that $v_i$ is a leaf in $\overline{G}_s^v$ for $2\leq i\leq s$. Set $e=\{v_s,c\}$. Then, by Corollary \ref{reg BEI for betti splitting} , and the fact that $(\overline{G}_s^v\setminus{e})_e=G_s^v$, we get 
   \begin{equation}\label{eq-star-clique-1}
       \reg(S/J_{\overline{G}_s^v})=\max\left\{\reg(S/J_{\overline{G}_{s-1}^v}),\reg(S/J_{G_s^v})+1\right\}.
   \end{equation} By induction, we have $\reg(S/J_{\overline{G}_{s-1}^v})=\reg(S/J_{G_{s-1}^v})+1$. Since $G_{s-1}^v$ is an induced subgraph of $G_s^v$, the result follows from Corollary \ref{reg reln b/w inducedsubgraph} and Equation \eqref{eq-star-clique-1}.
\end{proof}

Let $G$ be a graph, and $v\in V(G)$ a leaf vertex with the adjacent vertex $w$. For natural numbers $s\geq 1$ and $t\geq 2$, denote $G_{s,t}^v$ the graph obtained by identifying $v$ with one vertex of $K_t$, and attaching whiskers $e_1,\ldots,e_s$ at $v$, i.e., $G_{s,t}^v=G\cup _vK_t\cup_v\underbrace{K_2\cup_v\cdots \cup_vK_2}_{s\text{ times}}$.

The following figure illustrates the construction of $G_{s,t}^v.$
\begin{center}
\begin{tikzpicture}[scale=.35, line cap=round,line join=round,x=1.0cm,y=1.0cm]
\clip(-10.585681940000006,-4.946980480000011) rectangle (20.4921661000000075,4.358366320000012);
\draw (0,0)-- (0,2);
\draw (0,0)-- (-2,0);
\draw (0,0)-- (0,-2);
\draw (0,0)-- (2,0);
\draw (2,0)-- (0.81,1.56);
\draw (2,0)-- (3.18,1.56);
\draw (0.81,1.56)-- (1.95,3.18);
\draw (1.95,3.18)-- (3.18,1.56);
\draw (1.95,3.18)-- (2,0);
\draw (0.81,1.56)-- (3.18,1.56);
\draw (2,0)-- (4,0);
\draw (2,0)-- (2,-2);

\draw (-8,0)-- (-6,0);
\draw (-8,0)-- (-8,2);
\draw (-8,0)-- (-10,0);
\draw (-8,0)-- (-8,-2);
 \draw (12,0)-- (12,2);
 \draw (12,0)-- (12,-2);
\draw (10,0)-- (12,0);
\draw (12.78,1.98)-- (12,0);
\draw (12.78,1.98)-- (14.82,2.01);
\draw (14.82,2.01)-- (16,0);
\draw (16,0)-- (14.91,-1.95);
\draw (14.91,-1.95)-- (12.9,-2.01);
\draw (12.9,-2.01)-- (12,0);
\draw (12.78,1.98)-- (16,0);
\draw (12.78,1.98)-- (14.91,-1.95);
\draw (12.78,1.98)-- (12.9,-2.01);
\draw (14.82,2.01)-- (12,0);
\draw (14.82,2.01)-- (12.9,-2.01);
\draw (14.82,2.01)-- (14.91,-1.95);
\draw (16,0)-- (12,0);
\draw (16,0)-- (12.9,-2.01);
\draw (14.91,-1.95)-- (12,0);
\begin{scriptsize}
\fill (0,0) circle (3.5pt);
\fill (0,2) circle (3.5pt);
\fill (-2,0) circle (3.5pt);
\fill (0,-2) circle (3.5pt);
\fill (2,0) circle (3.5pt);
\fill (0.81,1.56) circle (3.5pt);
\fill (3.18,1.56) circle (3.5pt);
\fill (1.95,3.18) circle (3.5pt);
\fill (4,0) circle (3.5pt);
\fill (2,-2) circle (3.5pt);

\fill (-8,0) circle (3.5pt);
\fill (-6,0) circle (3.5pt);
\fill (-8,2) circle (3.5pt);
\fill (-10,0) circle (3.5pt);
\fill (-8,-2) circle (3.5pt);
\fill (10,0) circle (3.5pt);
\fill (12,2) circle (3.5pt);
\fill (12,-2) circle (3.5pt);
\fill (12,0) circle (3.5pt);
\fill (12.78,1.98) circle (3.5pt);
\fill (14.82,2.01) circle (3.5pt);
\fill (16,0) circle (3.5pt);
\fill (14.91,-1.95) circle (3.5pt);
\fill (12.9,-2.01) circle (3.5pt);
\end{scriptsize}
\draw[color=black] (-8.5,-.35) node {$w$};
\draw[color=black] (-6.4,-.35) node {$v$};
\draw[color=black] (-.5,-.35) node {$w$};
\draw[color=black] (1.6,-.35) node {$v$};
\draw[color=black] (11.5,-.35) node {$w$};
\draw[color=black] (-8,-4) node {$G$};
\draw[color=black] (1.2,-4) node {$G_{2,4}^v$};
\draw[color=black] (13.2,-4) node {$(G\setminus{v})_{6}^w$};
\end{tikzpicture}
\end{center}
\begin{lemma}\label{gluing of on clique and edges to a leaf}
Adopt the notation as above. Then, for all $s\geq 1$ and $t\geq 2$, we have $$\reg\left(S/J_{G_{s,t}^v}\right)=1+\reg\left(S/J_{(G\setminus{v})_{s+t}^w}\right).$$
\end{lemma}
\begin{proof}
Let $e_1,\ldots,e_s$ be edges and $e_s=\{i,j\}$. Without loss of generality, assume that $i=v$, and $j$ is a leaf vertex in $G_{s,t}^v$. Then the fact $(G_{s,t}^v\setminus{e_s})_{e_s}=(G\setminus{v})_{s+t}^w$ and Corollary \ref{reg BEI for betti splitting} give
    \begin{equation}\label{eq-edges-clique-1}
       \reg(S/J_{G_{s,t}^v})=\max\left\{\reg(S/J_{G_{s-1,t}^v}),\reg(S/J_{(G\setminus{v})_{s+t}^w})+1\right\}.
    \end{equation}
 In order to prove the lemma, we proceed by induction on $s$. For $s=1$, Equation \eqref{eq-edges-clique-1} gives
    $\reg(S/J_{G_{1,t}^v})=\max\{\reg(S/J_{G_t^v}),\reg(S/J_{(G\setminus{v})_{1+t}^w})+1\}$. 
    Note that $G_t^v=G\cup K_t$ is a decomposition of $G$. Then, using Theorem \ref{Gluing theorem}, we obtain $\reg(S/J_{G_t^v})=\reg(S/J_G)+1$. Since $G$ is an induced subgraph of $(G\setminus{v})_{1+t}^w$, the result follows from Corollary \ref{reg reln b/w inducedsubgraph}. 
    Now, assume $s\geq 2$.
        Then, by induction, we obtain $\reg(S/J_{G_{s-1,t}^v})=1+\reg(S/J_{(G\setminus{v})_{s-1+t}^w})$. Since the graph  $(G\setminus{v})_{s-1+t}^w$ is an induced subgraph of $(G\setminus{v})_{s+t}^w$,  by Corollary \ref{reg reln b/w inducedsubgraph}, we get $\reg(S/J_{(G\setminus{v})_{s-1+t}^w})\leq \reg(S/J_{(G\setminus{v})_{s+t}^w})$. Thus, the lemma follows from Equation \eqref{eq-edges-clique-1}. 
\end{proof}
The following example illustrates Lemma \ref{gluing Ks and star} and Lemma \ref{gluing of on clique and edges to a leaf}.
 \begin{example}\label{exm: caterpillar replaced by star}
     Let $G, G_1$ and $G_2$ be graphs shown in the Figure \ref{fig: for caterpillar attached to a vertex}. Then, by Lemma \ref{gluing Ks and star}, we have $\reg(S/J_{G})=1+\reg(S/J_{G_1})$. Further, applying Lemma \ref{gluing of on clique and edges to a leaf} and then Lemma \ref{gluing Ks and star}, we get $\reg(S/J_{G_1})=\reg(S/J_{G_2})$. Hence $\reg(S/J_G)=1+\reg(S/J_{G_2})$.
 \end{example}
  \begin{figure}[H]
     \centering
    \begin{tikzpicture}[scale=.35, line cap=round,line join=round,x=1.0cm,y=1.0cm]
\draw (0,0)-- (0,-2);
\draw (0,2)-- (0,0);
\draw (-2,0)-- (0,0);
\draw (2,0)-- (0,0);
\draw (-2,0)-- (-3.18,1.49);
\draw (-2,0)-- (-3.18,-1.49);
\draw (-4,0)-- (-5.18,-1.49);
\draw (-4,0)-- (-5.18,1.49);
\draw (-4,0)-- (-2,0);
\draw (2,0)-- (3.18,1.49);
\draw (2,0)-- (3.18,-1.49);
\draw (0,2)-- (-1.18,3.49);
\draw (0,2)-- (1.18,3.49);
\draw (0,2)-- (0,4);
\draw (0,4)-- (-1.18,5.49);
\draw (0,4)-- (1.18,5.49);
\draw[color=black] (0.35,-0.40) node {$u$};
\draw[color=black] (0,-3.60) node {$G$};
\draw (12,0)--(14,0);
\draw (12,0)--(10,0);
\draw (8,0)-- (10,0);
 \draw (10,0)--(8.82,1.49);
\draw (10,0)--(8.82,-1.49);
 \draw (14,0)--(15.18,1.49);
 \draw (14,0)--(15.18,-1.49);
 \draw (12,0)--(12,2);
 \draw (12,2)--(10.82,3.49);
\draw (12,2)--(13.18,3.49);
 \draw (12,2)--(12,4);
 \draw (8,0)-- (6.82,-1.49);
 \draw (8,0)-- (6.82,1.49);
  \draw (12,-2)-- (12,0);
  \draw (12,4)-- (10.82,5.59);
  \draw (12,2)-- (10.82,5.59);
 \draw[color=black] (12.35,-0.4) node {$u$};
 \draw[color=black] (12,-3.6) node {$G_1$};
 \draw[color=black] (0.4,1.6) node {$v$};
 \draw[color=black] (0.4,3.6) node {$w$};
  \draw[color=black] (12.4,1.6) node {$v$};
  \draw[color=black] (12.4,3.6) node {$w$};
  \draw[color=black] (12.4,3.6) node {$w$};
  \draw (24,0)--(26,0);
\draw (24,0)--(24,-2);
\draw (20,0)--(18.82,-1.49);
\draw (20,0)--(18.82,1.49);
\draw (24,0)--(22,0);
\draw (20,0)--(22,0);
\draw (22,0)--(20.82,1.49);
\draw (22,0)--(20.82,-1.49);
\draw (26,0)--(27.18,1.49);
\draw (26,0)--(27.18,-1.49);
\draw (24,0)--(24,2);
\draw (24,2)--(22.4,3);
\draw (24,2)--(25.6,3);
\draw (24,2)--(23.2,3.7);
\draw (24,2)--(24.8,3.7);
 \draw[color=black] (24.4,-.4) node {$u$};
 \draw[color=black] (24.4,-3.6) node {$G_2$};
 \begin{scriptsize}
 \fill (10.82,5.59) circle (3.5pt);
     \fill (12,-2) circle (3.5pt);
     \fill (0,0) circle (3.5pt);
\fill (0,-2) circle (3.5pt);
\fill (0,2) circle (3.5pt);
\fill (-2,0) circle (3.5pt);
\fill (2,0) circle (3.5pt);

\fill (-3.18,1.49) circle (3.5pt);
\fill (-3.18,-1.49) circle (3.5pt);

\fill (-4,0) circle (3.5pt);
\fill (-5.18,1.49) circle (3.5pt);
\fill (-5.18,-1.49) circle (3.5pt);

\fill (3.18,1.49) circle (3.5pt);
\fill (3.18,-1.49) circle (3.5pt);

\fill (-1.18,3.49) circle (3.5pt);
\fill (1.18,3.49) circle (3.5pt);

\fill (0,4) circle (3.5pt);
\fill (-1.18,5.49) circle (3.5pt);
\fill (1.18,5.49) circle (3.5pt);

\fill (12,0) circle (3.5pt);
\fill (14,0) circle (3.5pt);
\fill (10,0) circle (3.5pt);
\fill (8,0) circle (3.5pt);

\fill (8.82,1.49) circle (3.5pt);
\fill (8.82,-1.49) circle (3.5pt);

\fill (15.18,1.49) circle (3.5pt);
\fill (15.18,-1.49) circle (3.5pt);

\fill (12,2) circle (3.5pt);
\fill (10.82,3.49) circle (3.5pt);
\fill (13.18,3.49) circle (3.5pt);
\fill (12,4) circle (3.5pt);
\fill (6.82,-1.49) circle (3.5pt);
\fill (6.82,1.49) circle (3.5pt);

\fill (24,0) circle (3.5pt);
\fill (26,0) circle (3.5pt);
\fill (24,-2) circle (3.5pt);
\fill (20,0) circle (3.5pt);
\fill (18.82,-1.49) circle (3.5pt);
\fill (18.82,1.49) circle (3.5pt);
\fill (22,0) circle (3.5pt);
\fill (20.82,1.49) circle (3.5pt);
\fill (20.82,-1.49) circle (3.5pt);
\fill (27.18,1.49) circle (3.5pt);
\fill (27.18,-1.49) circle (3.5pt);
\fill (24,2) circle (3.5pt);
\fill (22.4,3) circle (3.5pt);
\fill (25.6,3) circle (3.5pt);
\fill (23.2,3.7) circle (3.5pt);
\fill (24.8,3.7) circle (3.5pt);
 \end{scriptsize}
\end{tikzpicture}
\caption{}\label{fig: for caterpillar attached to a vertex}
\end{figure}
 Recall that a tree is \emph{caterpillar} if the induced subgraph on internal vertices is either empty or a path. Example \ref{exm: caterpillar replaced by star} shows that the induced caterpillar on the vertex set $N_G[v]\cup N_G[w]$ can be replaced by a star graph while studying the regularity of $G$. More generally, we have the following:
\begin{lemma}\label{caterpillar attached to a vertex of a graph}
Let $G$ be a graph, and there are induced subgraphs $H$ and caterpillar $T$ with spine $P:v=v_l,v_{l-1},\ldots,v_1,v_0$ of length $l\geq 3$ such that $G=H\cup_vT$. If $\deg(v_i)\geq 3$ for $1\leq i\leq l-1$, then $$\reg(S/J_G)=l-2+\reg(S/J_{\overline{H}_{t}^v}),\text{ where } t=\sum\limits_{i=1}^{l-1}\deg(v_i)-2l+4.$$ 
\end{lemma}
\begin{proof}  Let $P:v=v_l,v_{l-1},\ldots,v_{2},v_1,v_0$ be a spine of length $l\geq 3$ of caterpillar $T$ such that $G=H\cup_vT$. Set $\deg(v_i)=d_i$, for $1\leq i\leq l-1$. Note that graph $G$ can be seen as $(\overline{G\setminus{v_1}})_{d_1}^{v_2}$. Therefore, by Lemma \ref{gluing Ks and star}, we have $\reg(S/J_G)=1+\reg(S/J_{(G\setminus{v_1})_{d_1}^{v_2}})$. One can see that $(G\setminus{v_1})_{d_1}^{v_2}$ satisfies the hypothesis of Lemma \ref{gluing of on clique and edges to a leaf} with $s=d_{2}-2$ and $t=d_{1}$. This implies $\reg(S/J_{(G\setminus{v_{1}})_{d_{1}}^{v_{2}}})=1+\reg(S/J_{(G\setminus{v_{1},v_{2}})_{s'}^{v_{3}}})$, where $s'={d_{1}+d_{2}-2}$. Using Lemma \ref{gluing Ks and star} for the graph $(G\setminus{v_{1},v_{2}})_{s'}^{v_{3}}$, we obtain $$\reg\left(S/J_{\overline{(G\setminus{v_{1},v_{2}})}_{s'}^{v_{3}}}\right)=1+\reg\left(S/J_{(G\setminus{v_{1},v_{2}})_{s'}^{v_{3}}}\right)=\reg(S/J_{(G\setminus{v_{1}})_{d_{1}}^{v_{2}}}).$$ Thus, we have $\reg(S/J_G)=1+\reg(S/J_{{\overline{(G\setminus{v_{1},v_{2}})}_{s'}^{v_{3}}}})$. We use induction on $l$ to prove the lemma.  For $l=3$, observe that $\overline{(G\setminus{v_{1},v_{2}})}_{s'}^{v_{3}}$ is the graph $\overline{H}_{s'}^v$, and hence the result follows. Now, assume $l\geq 4$. Then, the graph $\overline{(G\setminus{v_{1},v_{2})}}_{s'}^{v_{3}}=H\cup_vT'$ satisfies the hypothesis of the lemma with caterpillar $T'$ that has a spine $P':v=v_l,v_{l-1},\ldots,v_3,v_2',v_1'$ of length $l-1$. Therefore, by induction, we have $$\reg(S/J_{\overline{(G\setminus{v_{1},v_{2}})}_{s'}^{v_{3}}})=(l-3)+\reg(S/J_{\overline{H}_{t}^v}),$$ where $t=\deg(v_2')+\sum\limits_{i=3}^{l-1}\deg(v_i)-2(l-1)+4$. Since $\deg(v_2')=d_1+d_2-2$, we obtain $t=\sum\limits_{i=1}^{l-1}\deg(v_i)-2l+4$. Thus, we get $$\reg(S/J_G)=1+\reg(S/J_{{\overline{(G\setminus{v_{1},v_{2}})}_{s'}^{v_{3}}}})=l-2+\reg(S/J_{\overline{H}_t^v}).$$ This completes the proof.
\end{proof}
As an immediate consequence of Lemma \ref{caterpillar attached to a vertex of a graph}, and using the fact that $\reg(S/J_{K_{1,r}})=2$ for $r\geq 2$, we get the following, which was proved in \cite[Theorem 4.1]{On-the-BEI-of-Block-graphs-Chaudhary-Irfan}.
\begin{corollary}
    For a caterpillar $G$ of spine length $l$, we have $\reg(S/J_G)=l$.
\end{corollary}

The result outlined below is a slightly strengthened version of \cite[Theorem 4.4]{RegBEIcertainblockgraph}. We demonstrate it through the application of Betti splitting, which yields a similar proof. We recall that $\mathcal{G}(s,t,m;c)$ (cf. \cite{RegBEIcertainblockgraph}) is the family of graphs obtained by identifying a leaf vertex of each of $K_{1,t_1},\ldots, K_{1,t_s}$, where $t_i\geq 3,1\leq i \leq s$, $t$ cliques on at least three vertices and $m$ whiskers at $c$.
\begin{lemma}\label{reg of stars at a point}
    Let $G\in \mathcal{G}(s,t,m;c)$. Then, for all $s+t\geq 2$, we have $$\reg\left(S/J_{G}\right)=2s+t.$$

\end{lemma}
\begin{proof}
   We prove the lemma by induction on $s$. If $s=0$, then $G\in \mathcal{G}(0,t,m;c)$ with $t\geq 2$. By Corollary \ref{reg reln b/w inducedsubgraph}, we get $\reg(S/J_{G})\geq t$ . Let $G', G''$ and $\Tilde{G}$ be the graphs obtained by applying the EHH-process on $G$ with respect to $c$. Then, $\reg(S/J_{G'})=\reg(S/J_{\Tilde{G}})=1$. Since $G''$ is a disjoint union of $t$ cliques and $m$ isolated vertices, we have $\reg(S/J_{G''})=t$. Therefore, the assertion follows from the short exact sequence:
   $$0\rightarrow{S/J_G}\rightarrow S/J_{G'}\oplus S/J_{G''}\rightarrow S/J_{\Tilde{G}}\rightarrow 0.$$
   Assume $s\geq 1$. By Lemma \ref{gluing Ks and star}, we have $$\reg(S/J_{G})=1+\reg(S/J_{G_1}),$$ 
   where $G_1\in \mathcal{G}(s-1,t+1,m;c)$. 
   Note that $G_1$ satisfies the hypothesis of the lemma with $s-1$ stars attached to $c$. By induction, we have $\reg(S/J_{G_1})=2(s-1)+(t+1)$. Thus, we obtain $\reg(S/J_{G})=1+2(s-1)+t+1=2s+t$.
\end{proof}

\section{regularity of trees}\label{SecMainResults}
   This section studies the regularity of binomial edge ideals for trees. Let $G$ be a connected tree on $[n]$ and $J_G$ its binomial edge ideal in $S$. Recall that the jewel $\mathcal J $ is a tree obtained by identifying a leaf vertex of each of three copies of $K_{1,3}$, and that vertex is called the center of $\J$. Jayanthan et al. \cite[Theorem 4.1, Theorem 4.2]{RegBEIcertainblockgraph} proved that $\reg(S/J_G)\geq\iv(G)+1$ with equality if and only if $G$ is jewel-free.
   
Let $G$ be a tree that contains a jewel. Then one can observe that $G$ contains an induced subgraph $H\in \mathcal G(s,0,m;c)$ for some $s\geq 3$ and $m\geq 0$. The regularity of this family of graphs is computed in Lemma \ref{reg of stars at a point}. We observe that, within this family, the regularity depends only on the parameter $s$. In particular, for every graph $H\in \mathcal G(3,0,m;c)$, we have $\reg(S/J_H)=\reg(S/J_{\mathcal J})$. The induced jewel $\J$ of $G$ with the same center may not be unique. It turns out that a maximal induced subgraph of $G$ contained in $\cup_{s\geq 3} \mathcal{G}(s,0,0;c)$ is unique, and it contains all the jewels with the same center. The maximal induced subgraph of $G$ contained in $\cup_{s\geq 3} \mathcal{G}(s,0,0;c)$ is called a generalized jewel. We give a more formal definition of a generalized jewel as follows:

For $v\in V(G)$ and $i\in \mathbb N$, denote $N_{G}^{\geq i}(v)$ the set of adjacent vertices to $v$ of degree $\geq i$. Set $D_G({v}):=|N_{G}^{\geq 3}(v)|$. When $D_G(v)\geq 3$,  we call the induced subgraph on $\{N_G[u]: u\in N_{G}^{\geq 3}(v)\}$ a \emph{generalized jewel} or simply \emph{jewel}, denoted $\mathcal{J}_v$. The vertex $v$ is the \emph{center}, and $u\in N_G^{\geq 3}(v)$ is a \emph{supporting vertex} of the jewel $\mathcal{J}_v$. Note that all the supporting vertices are internal vertices. For simplicity of notation, we omit subscripts and write $D(v)$ and $N^{\geq i}(v)$ instead of $D_G(v)$ and $N_G^{\geq i}(v)$ when it is clear from the context. Denote $D_G = \sum D_G(c)$, where the summation runs over the centers of jewels.

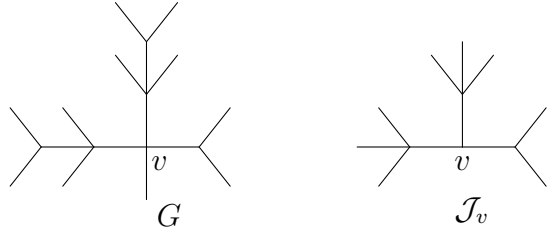
\begin{figure}[H]
     \centering
    \begin{tikzpicture}[scale=.35, line cap=round,line join=round,x=1.0cm,y=1.0cm]
\draw (0,0)-- (0,-2);
\draw (0,2)-- (0,0);
\draw (-2,0)-- (0,0);
\draw (2,0)-- (0,0);
\draw (-2,0)-- (-3.18,1.49);
\draw (-2,0)-- (-3.18,-1.49);
\draw (-4,0)-- (-5.18,-1.49);
\draw (-4,0)-- (-5.18,1.49);
\draw (-4,0)-- (-2,0);
\draw (2,0)-- (3.18,1.49);
\draw (2,0)-- (3.18,-1.49);
\draw (0,2)-- (-1.18,3.49);
\draw (0,2)-- (1.18,3.49);
\draw (0,2)-- (0,4);
\draw (0,4)-- (-1.18,5.49);
\draw (0,4)-- (1.18,5.49);
\draw[color=black] (0.50,-0.50) node {$v$};
\draw[color=black] (0,-4.60) node {$G$};
\draw (12,0)--(14,0);
\draw (12,0)--(10,0);
\draw (8,0)-- (10,0);
 \draw (10,0)--(8.82,1.49);
\draw (10,0)--(8.82,-1.49);
 \draw (14,0)--(15.18,1.49);
 \draw (14,0)--(15.18,-1.49);
 \draw (12,0)--(12,2);
 \draw (12,2)--(10.82,3.49);
\draw (12,2)--(13.18,3.49);
 \draw (12,2)--(12,4);
 \draw (0,-2)-- (-1.18,-3.49);
 \draw (0,-2)-- (1.18,-3.49);
 \draw (12,-2)-- (10.82,-3.49);
  \draw (12,-2)-- (12,0);
  \draw (12,-2)-- (13.18,-3.49);
 \draw[color=black] (12.5,-0.5) node {$v$};
 \draw[color=black] (12,-4.6) node {$\mathcal{J}_v$};
 \begin{scriptsize}
 \fill (10.82,-3.49) circle (3.5pt);
     \fill (13.18,-3.49) circle (3.5pt);
     \fill (12,-2) circle (3.5pt);
     \fill (0,0) circle (3.5pt);
\fill (0,-2) circle (3.5pt);
\fill (0,2) circle (3.5pt);
\fill (-2,0) circle (3.5pt);
\fill (2,0) circle (3.5pt);

\fill (-3.18,1.49) circle (3.5pt);
\fill (-3.18,-1.49) circle (3.5pt);

\fill (-4,0) circle (3.5pt);
\fill (-5.18,1.49) circle (3.5pt);
\fill (-5.18,-1.49) circle (3.5pt);

\fill (3.18,1.49) circle (3.5pt);
\fill (3.18,-1.49) circle (3.5pt);

\fill (-1.18,3.49) circle (3.5pt);
\fill (1.18,3.49) circle (3.5pt);

\fill (0,4) circle (3.5pt);
\fill (-1.18,5.49) circle (3.5pt);
\fill (1.18,5.49) circle (3.5pt);

\fill (12,0) circle (3.5pt);
\fill (14,0) circle (3.5pt);
\fill (10,0) circle (3.5pt);
\fill (8,0) circle (3.5pt);

\fill (8.82,1.49) circle (3.5pt);
\fill (8.82,-1.49) circle (3.5pt);

\fill (15.18,1.49) circle (3.5pt);
\fill (15.18,-1.49) circle (3.5pt);

\fill (12,2) circle (3.5pt);
\fill (10.82,3.49) circle (3.5pt);
\fill (13.18,3.49) circle (3.5pt);

\fill (12,4) circle (3.5pt);
\fill (-1.18,-3.49) circle (3.5pt);
\fill (1.18,-3.49) circle (3.5pt);
 \end{scriptsize}
\end{tikzpicture}
\caption{A tree $G$ containing a generalized jewel $\mathcal J_v$}\label{fig: generalized jewel}
\end{figure}

The tree $G$ in the Figure \ref{fig: generalized jewel} has only one vertex, namely $v$, with $D(v)\geq3$. More precisely $D(v)=4$. Observe that $G$ contains more than one jewel $\mathcal{J}$ as an induced subgraph of $G$. In fact, it contains induced subgraphs from the family $\mathcal G(s,0,0;v)$ for $s=3,4$. The maximal induced subgraph of $G$ in $\mathcal G(3,0,0;v)\cup\mathcal G(4,0,0;v)$ is $\mathcal{ J}_v$. By definition, this is a generalized jewel with center $v$. 

The following remark contains some basic observations that will be useful throughout the section.
\begin{remark}\label{reln-split-dG-iv-jewels}\hfill {}
\begin{enumerate}
    \item \label{one jewel degree condition} Let $G$ be a tree containing exactly one generalized jewel $\mathcal{J}_c$ centered at $c$. Then $D(v)\leq 2$ for all $v\in V(G)\setminus{\{c\}}$.
    \item \label{split reln of iv and D-value}Let $v\in V(G)$ with $\deg(v)=2$. Then the decomposition $G=G_1\cup G_2$ of $G$ at $v$ gives that $\iv(G)=\iv(G_1)+\iv(G_2)+1$. If a vertex $u(\neq v)\in V(G_i)$, then $N_G^{\geq 3}(u)\subset V(G_i)$. Therefore, a jewel $\mathcal{J}_u$ is contained in either $G_1$ or $G_2$, and $D_G(u)=D_{G_i}(u)$ for $G_i$ that contains $u$. Consequently, if $s_i$ is the number of jewels in $G_i$ for $i=1,2$, then $s_1+s_2=s$, and $D_G=D_{G_1}+D_{G_2}$, where $s$ is the number of jewels in $G$.
\end{enumerate}
\end{remark}
In the lemma below, we establish the regularity of trees containing precisely one jewel.
\begin{lemma}\label{reg-one-jewel}
        Let $G$ be a tree containing exactly one generalized jewel $\mathcal{J}_c$ centered at $c$. Then $$\reg(S/J_G)=\iv(G)+D(c)-1.$$
\end{lemma}
\begin{proof}
Suppose there exists a vertex $v\in V(G)$ such that $\deg(v)=2$. Then $G=G_1\cup G_2$ is the decomposition of $G$ at $v$. By Theorem \ref{Gluing theorem}, we get $\reg(S/J_G)=\reg(S/J_{G_1})+\reg(S/J_{G_2})$. In view of Remark \ref{reln-split-dG-iv-jewels}\eqref{split reln of iv and D-value}, assume that $\mathcal{J}_c$ is contained in $G_1$. Then the graph $G_2$ is jewel-free as $G$ contains exactly one generalized jewel, and $\reg(S/J_{G_2})=\iv(G_2)+1$. The assertion follows from Remark \ref{reln-split-dG-iv-jewels}\eqref{split reln of iv and D-value} if it holds for $G_1$. 

 Now, we can assume that $G$ does not have a vertex of degree two. Since $c$ is an internal vertex, we have $\iv(G)\geq D(c)+1$. We prove the lemma by induction on $\iv(G)$.  For $\iv(G)=D(c)+1$, note that $G\in \mathcal G(D(c),0,m;c)$ for $m\geq 0$, and the result follows from Lemma \ref{reg of stars at a point}.
We now assume that $ \iv(G)>D(c)+1$. Set $l:=\max\{\dist(c,u):\text{$u$ is a leaf in $G$}\}$. Let $v$ be a leaf such that $l=\dist(c,v)$. Since $ \iv(G)>D(c)+1$, and there is no vertex of degree two in $G$, we have $l\geq 3$. Let $P:c=v_l,v_{l-1},\ldots,v_1,v$ be the path from $c$ to $v$, and $\deg(v_i)=d_i$ for $i=1,\ldots,l-1$. Since $d_i\geq 3$ and $v_i\neq c$ for $i=1,\ldots,l-1$, all the adjacent vertices of $v_i$ except $v_{i-1}$ and $v_{i+1}$ are leaves by Remark \ref{reln-split-dG-iv-jewels}\eqref{one jewel degree condition}. Denote $T$ the caterpillar induced on $\bigcup\limits_{i=1}^{l-1}N[v_i]$, and $H$ the induced subgraph of $G$ on $(V(G)\setminus{V(T)})\cup \{c\}$. Then $G=H\cup_{c}T$, and by Lemma \ref{caterpillar attached to a vertex of a graph}, we get \begin{equation}\label{eq 6}
\reg(S/J_G)=l-2+\reg(S/J_{\overline{H}_{t}^{c}})\end{equation} with $t=\sum\limits_{i=1}^{l-1}d_i-2l+4$. Notice that $ t\geq 4$. Since $\overline{H}_{t}^{c}$ is obtained by identifying a leaf vertex of $K_{1,t}$ at $c$, the center of $K_{1,t}$ belongs to $N^{\geq 3}(c)$ in $\overline{H}_{t}^{c}$. Other supporting vertices of $\mathcal{J}_c$ remains unaltered in $\overline{H}_{t}^{c}$. Therefore, one has $D_{\overline{H}_{t}^{c}}(c)=D_G(c)$, and the graph $\overline{H}_{t}^{c}$ is a tree containing exactly one jewel with $\iv(\overline{H}_{t}^{c})=\iv(G)-(l-2)<\iv(G)$. 

By induction, we obtain $$\reg\left(S/J_{\overline{H}_{t}^c}\right)=\iv(G)-(l-2)+D(c)-1.$$ Thus, the desired result follows from Equation \eqref{eq 6}. 
    \end{proof}
Recall that a spine of a tree $G$ is a path of maximum length in $G.$ The following remark ensures the location of the center of some jewel.
\begin{remark}\label{existance-of-jewel-center-on spine}
    Let $G$ be a tree containing a jewel. Suppose that there is no vertex of degree two in $G$, and that for every caterpillar $T$ in $G$ with spine $w=w_l,\ldots,w_0$ such that $G=H\cup_w T$ for some $H$, the spine length of $T$ is at most two. Then, for every spine $P:v_l,\ldots,v_0$ of $G$, the vertices $v_2$ and $v_{l-2}$ are centers of jewels. 
\end{remark}
We now recall the previously introduced notation and introduce additional technical notation that will be used in the subsequent results.

\begin{notation}
    Let $G$ be a tree containing a jewel. Recall that for $v\in V(G)$, $D(v)$ denotes the number of adjacent vertices to $v$ having degree at least three. A vertex $v$ is the center of a jewel if $D(v)\geq 3$,  and the number $D_G=\sum D(v)$, where summation runs over the center of jewels.  

    Let $G_{\mathcal{J}}$ denote the induced subgraph of $G$ on the vertex set $\{c:c\text{ is a center of a jewel in $G$}\}$. Denote $G_1,\ldots, G_p$ the connected components of $G_{\mathcal{J}}$, and $s_i=|V(G_i)|$ for $1\leq i\leq p$. Then, we have $s=\sum\limits_{i=1}^{p}s_i$. The number of isolated edges in $G_{\J}$ is denoted by $e_G$. Further, define $$C_G:=\{v\in V(G):\deg(v)=3, D(v)=2, \text{ and } N^{\geq 3}(v)\text{ consists of centers of jewels}\}.$$ Denote $\mu_G$ the number of elements in $C_G$, i.e., $\mu_G=|C_G|$. It is easy to observe that if $G=G_1\cup G_2$ is a decomposition of $G$, and $\mu_i=|C_{G_i}|$ for $i=1,2$, then $\mu=\mu_1+\mu_2$.

Let $G$ be the graph shown below. Observe that $u,v$, and $z$ are centers of jewels in $G$. Consequently, $D_G=12$. The graph $G_{\J}$ consists of an edge and an isolated vertex. More precisely, $G_{\J}=\{\{u,v\},\{z\}\}$, and $C_G=\{w\}$. Therefore, we have $p=2$, $e_G=1$, and $\mu_G=1$. 
\end{notation}
\begin{figure}[ht]
\centering

\begin{tikzpicture}[scale=0.45,line cap=round,line join=round,x=1cm,y=1cm]

\draw (0,0)--(0,2);
\draw (0,0)--(0,-2);
\draw (0,0)--(2,0);
\draw (2,0)--(2,2);
\draw (2,0)--(2,-2);
\draw (0,0)--(-2,0);
\draw (2,0)--(4,0);

\draw (0,2)--(-0.6,3);
\draw (0,2)--(0.6,3);
\draw (2,2)--(1.4,3);
\draw (2,2)--(2.6,3);

\draw (-2,0)--(-3,0.6);
\draw (-2,0)--(-3,-0.6);

\draw (0,-2)--(-0.6,-3);
\draw (0,-2)--(0.6,-3);

\draw (2,-2)--(1.4,-3);
\draw (2,-2)--(2.6,-3);
\draw (4,0)--(4,2);
\draw (4,0)--(6,0);
\draw (6,0)--(6,2);
\draw (6,0)--(6,-2);
\draw (6,0)--(8,0);

\draw (6,2)--(5.4,3);
\draw (6,2)--(6.6,3);

\draw (6,-2)--(5.4,-3);
\draw (6,-2)--(6.6,-3);

\draw (8,0)--(9,0.6);
\draw (8,0)--(9,-0.6);

\foreach \x/\y in {
0/0,0/2,0/-2,
2/0,2/2,2/-2,
-2/0,
4/0,
4/2,
6/0,6/2,6/-2,
8/0,
-0.6/3,0.6/3,
1.4/3,2.6/3,
-3/0.6,-3/-0.6,
-0.6/-3,0.6/-3,
1.4/-3,2.6/-3,
5.4/3,6.6/3,
5.4/-3,6.6/-3,
9/0.6,9/-0.6}
{
\fill (\x,\y) circle (3.5pt);
}

\node at (-0.4,-0.4) {$u$};
\node at (1.5,-0.4) {$v$};
\node at (3.6,-0.4) {$w$};
\node at (6.4,-0.4) {$z$};
\end{tikzpicture}

\caption{}
\label{fig:merged}
\end{figure}
We now provide an upper bound for the regularity of binomial edge ideals of trees.
\begin{theorem}\label{upper bound for reg}
    Let $G$ be a tree containing $s\geq 1$ generalized jewels. Suppose $G_1,\ldots,G_p$ are connected components of $G_{\mathcal{J}}$, and $|V(G_i)|=s_i$ for all $1\leq i\leq p$. Denote $e_G$ the number of isolated edges in $G_{\mathcal{J}}$. Then
        $$\reg(S/J_G)\leq \iv(G)+1+D_G-2s-e_G-\sum\limits_{i=1}^{p}\left\lfloor \frac{s_i}{3}\right\rfloor.$$
\end{theorem}
\begin{proof}
    Using Theorem \ref{Gluing theorem} and Remark \ref{reln-split-dG-iv-jewels}\eqref{split reln of iv and D-value}, we assume that $G$ does not have a vertex of degree two. In view of Lemma \ref{caterpillar attached to a vertex of a graph}, assume that every caterpillar $T$ with spine $w=w_l,\ldots,w_0$ such that $G=H\cup_w T$ for some $H$, the spine length of $T$ is at most two.

 We prove the assertion using induction on $s$. The case $s=1$ follows from Lemma \ref{reg-one-jewel}. Now, assume $s\geq 2$. Then, any spine of $G$ has length at least $5$. Let $P:v_l,v_{l-1},\ldots,v_0$ be a spine in $G$. Then, by Remark \ref{existance-of-jewel-center-on spine}, $v_2$ is the center of a jewel. Note that all supporting vertices $v_1,w_3,\ldots,w_{D(v_2)}$ of $\mathcal{J}_{v_2}$ except $v_3$ are centers of star graphs on at least four vertices attached to ${v_2}$. Lemma \ref{gluing Ks and star} applied to these $D({v_2})-1$ stars attached to $v_2$ gives 
 \begin{equation}\label{proof-UB-stars-replaced-by-cliques}
     \reg(S/J_G)=D(v_2)-1+\reg(S/J_{G_{\Delta}}),
 \end{equation} where $G_{\Delta}$ is the graph obtained by replacing the star graphs attached to $v_2$ by cliques on $\deg(v_1),\deg(w_3),\ldots,\deg(w_{D(v_2)})$ vertices, respectively. Let $G_{\Delta}', G_{\Delta}''$ and $\Tilde{G_{\Delta}}$ be the graphs obtained by applying the EHH-process on $G_{\Delta}$ with respect to $v_2$. Denote $H$ the non-isolated component of the induced subgraph of $G$ on $(V(G)\setminus{N[v_2]})\cup \{v_3\}$.
 Then $G_{\Delta}'=H_t^{v_3}$ and $\Tilde{G_{\Delta}}=H_{t-1}^{v_3}$, where $ t=|N(v_2)|-|N^{\geq 2}(v_2)|+\deg(v_1)+1+\sum\limits_{j\geq 3}^{D(v_2)}(\deg(w_j)-1)$. Notice that $t\geq 6$. By Lemma \ref{gluing of on clique and edges to a leaf}, we have $$\reg(S/J_{G_{\Delta}'})=\reg(S/J_{\overline H_{t}^{v_3}})-1 \text{, and } \reg(S/J_{\Tilde{G_{\Delta}}})=\reg(S/J_{\overline H_{t-1}^{v_3}})-1.$$ The graph $G_{\Delta}''=G_{\Delta}\setminus{\{v_2\}}$ is a collection of $D(v_2)-1$ cliques on at least three vertices, $H$, and $\deg(v_2)-D(v_2)$ isolated vertices. Therefore, $\reg(S/J_{G_{\Delta}''})=D(v_2)-1+\reg(S/J_{H})$. The short exact sequence
		$$0\rightarrow S/J_{G_{\Delta}}\rightarrow S/J_{G_{\Delta}'}\oplus S/J_{G_{\Delta}''}\rightarrow S/J_{\Tilde{G_{\Delta}}}\rightarrow 0.$$ implies that $$
\reg(S/J_{G_{\Delta}})\leq \max\left\{\reg(S/J_{\overline{H}_t^{v_3}})-1,D(v_2)-1+\reg(S/J_H),\reg(S/J_{\overline{H}_{t-1}^{v_3}})\right\}
.$$
 Since $t\geq 6$, we have $|N_{\overline H_t^{v_3}}^{\geq 3}(v_3)|=|N_G^{\geq 3}(v_3)|$. Then $\overline H_{t}^{v_3}$ and $\overline H_{t-1}^{v_3}$ are trees containing $s-1$ jewels, and all the required numerical invariants are the same for these two graphs. Therefore, by induction and the above inequality, we have
  \begin{equation}\label{eqinequality G_delta}
      \reg(S/J_{G_{\Delta}})\leq \max\left\{D(v_2)-1+\reg(S/J_H),\reg(S/J_{\overline{H}_{t-1}^{v_3}})\right\}
  \end{equation}
  Since $v_2$ is the center of a jewel, there exists $i$ such that $v_2\in V(G_i)$, say $v_2\in V(G_p)$.
Now, we give an upper bound for $\reg(S/J_{\overline{H}_{t-1}^{v_3}})$. 
 
 Notice that $G_1,G_2,\ldots,G_{p-1},G_p'$ are the connected components of $(\overline{H}_{t-1}^{v_3})_{\mathcal J}$, where $G_p'=G_p\setminus{\{v_2\}}$, and we have $e_{\overline{H}_{t-1}^{v_3}}=e_G-1$ if $s_p=2$, and $e_{\overline{H}_{t-1}^{v_3}}\geq e_G$ if $s_p\neq 2$.
Since $\overline{H}_{t-1}^{v_3}$ contains $s-1$ jewels with $\iv(\overline{H}_{t-1}^{v_3})=\iv(G)-D({v_2})+1$ and $D_{\overline{H}_{t-1}^{v_3}}=D_G-D(v_2)$, by induction, we get
\begin{align*}
    \reg(S/J_{\overline{H}_{t-1}^{v_3}})\leq \lambda-2D(v_2)+3+\left\lfloor \frac{s_p}{3}\right\rfloor-\left\lfloor \frac{s_p-1}{3}\right\rfloor-
    \begin{cases}
        e_G-1, &\text{if }s_p=2,\\
        e_G, &\text{if } s_p\neq 2,
    \end{cases}
\end{align*}where $\lambda=\iv(G)+1+D_G-2s-\sum\limits_{i=1}^{p}\left\lfloor \frac{s_i}{3}\right\rfloor$. Since $\left\lfloor \frac{s_p-1}{3}\right\rfloor=0=\left\lfloor \frac{s_p}{3}\right\rfloor$ for $s_p=2$, and $\left\lfloor \frac{s_p-1}{3}\right\rfloor\geq\left\lfloor \frac{s_p}{3}\right\rfloor-1$ for $s_p\geq 3$, we obtain  
\begin{equation}\label{eq-bound for H_{t-1}}
 \reg(S/J_{\overline H_{t-1}^{v_3}})\leq\lambda-2D(v_2)+4-e_G.\end{equation}
 Next, we find an upper bound on $\reg(S/J_{H})$.  By applying Corollary \ref{reg reln b/w inducedsubgraph}, we assume that $\deg_H(v_3)\geq 3$, which allows us to derive the necessary bound on $\reg(S/J_G)$.
 Observe that $\iv(H)=\iv(G)-D({v_2})$. 
 
 Suppose $s_p=1$. Then $H$ contains $s-1$ jewels, and $D_{H}=D_G-D(v_2)$. The connected components of $H_{\mathcal J}$ are $G_1,\ldots,G_{p-1}$. Since $v_2$ is an isolated vertex in $G_{\mathcal J}$, we have $e_H=e_G$. Thus, by induction, we obtain $$\reg(S/J_{H})\leq \lambda-2D({v_2})+2-e_G.$$ Now, assume $s_p\geq 2$. In this case, $v_3$ is the center of a jewel. The following table summarizes the invariants of $H$ depending on the values of $D(v_3)$.
  \begin{table}[h]
\centering
\begin{tabular}{|c|c|c|c|c|}
\hline
 & no. of centers & $D_H$ & no. of components of $H_{\J}$ & $e_H$ \\ \hline
$D(v_3)\geq 4$& $s-1$& $D_G-D(v_2)-1$ & $p$& $e_G-1$ or $e_G$\\ \hline
$D(v_3)=3$ & $s-2$ & $D_G-D(v_2)-3$ & $p-1$, $p$ or $p+1$ & $e_G-1$ or $\geq e_G$ \\ \hline
\end{tabular}
\caption{}
\label{tab:example}
\end{table}
\newline
$\textbf{Case 1}:$ Suppose $D(v_3)\geq4$. Then $D_H(v_3)=D_G(v_3)-1\geq 3$. This means $v_3$ remains the center of a jewel in $H$. Therefore, $H$ contains $s-1$ jewels, and $D_H=D_G-D({v_2})-1$. The connected components of $H_{\mathcal J}$ are $G_1,\ldots,G_{p-1},G_p'$, where $G_p'=G_p\setminus \{v_2\}$, and $e_H=e_{\overline{H}_{t-1}^{v_3}}$. The same argument to obtain Equation \eqref{eq-bound for H_{t-1}} and induction gives the following: $$\reg(S/J_H)\leq \lambda-2D(v_2)+2-e_G.$$
$\textbf{Case 2}:$ Assume now that $D(v_3)=3$. Since $\deg_H(v_3)\geq 3$, the graph $H$ contains $s-2$ jewels, and $D_{H}=D_G-D({v_2})-3$. Let $w,v_4,v_2$ be supporting vertices of $\mathcal{J}_{v_3}$ in $G$. Then $G_1,\ldots G_{p-1},G_p',G_p''$ are the connected components of $H_{\mathcal J}$, where $G_p'=G_p\setminus{\{v_2,v_3,w\}}$ and $G_p''=\{w\}\cap V(G_p)$. If $w\in V(G_p)$, then $e_H\geq e_G$. Therefore, using $\left\lfloor \frac{s_p-3}{3}\right\rfloor=\left\lfloor \frac{s_p}{3}\right\rfloor-1$, we have $\reg(S/J_H)\leq \lambda-2D(v_2)+2-e_G$. Now, assume $w\notin V(G_p)$. In this case, $e_H= e_G-1$ if $s_p= 2$, and $e_H\geq e_G$ if $s_p\neq2$.
Then, by induction, we obtain $$\reg(S/J_H)\leq \lambda-2D(v_2)+1+\left\lfloor \frac{s_p}{3}\right\rfloor-\left\lfloor \frac{s_p-2}{3}\right\rfloor-e_H.$$
Considering all the cases for $s_p$, we have 
\begin{equation}\label{bound-J_H-in upper bound}
    \reg(S/J_H)\leq \lambda-2D(v_2)+2-e_G.
\end{equation}
By Equations \eqref{proof-UB-stars-replaced-by-cliques}, \eqref{eqinequality G_delta}, \eqref{eq-bound for H_{t-1}} and \eqref{bound-J_H-in upper bound}, we get 
$\reg(S/J_G)\leq \max\left \{ \lambda-e_G, \lambda-e_G-D(v_2)+3\right \}$. Since $D(v_2)\geq 3$, the assertion follows.
\end{proof}
We illustrate the above theorem with an example.
\begin{example}\label{ex: upper bound}
  In the graph $G$ shown below, we have $\iv(G)=8$, and there are two centers of jewels, namely $ v$ and $ w$. Therefore, we get $s=2$ and $D_G=8$. Note that $\{v,w\}$ is the only edge of the induced subgraph $G_{\mathcal{J}}$ on centers of jewels. Consequently, we get $p=1$, $s_1=s$, and $e_G=1$. Applying Theorem~\ref{upper bound for reg}, we conclude that $\reg(S/J_G)\leq 12$.  
\end{example} 
\begin{figure}[h]
\centering
\begin{tikzpicture}[scale=0.45, line cap=round,line join=round,,x=1cm,y=1cm]
\draw  (0,0)-- (0,2);
\draw  (0,0)-- (0,-2);
\draw  (0,0)-- (2,0);
\draw  (2,0)-- (2,2);
\draw  (2,0)-- (2,-2);
\draw  (0,0)-- (-2,0);
\draw  (2,0)-- (4,0);
\draw  (0,2)-- (-0.6,3);
\draw  (0,2)-- (0.6,3);
\draw  (2,2)-- (1.4,3);
\draw  (2,2)-- (2.6,3);
\draw  (-2,0)-- (-3,0.6);
\draw  (-2,0)-- (-3,-0.6);
\draw  (0,-2)-- (-0.6,-3);
\draw  (0,-2)-- (0.6,-3);
\draw  (2,-2)-- (1.4,-3);
\draw  (2,-2)-- (2.6,-3);
\draw  (4,0)-- (5,0.6);
\draw  (4,0)-- (5,-0.6);
\fill (0,0) circle (3.5pt);
\fill (2,0) circle (3.5pt);
\begin{scriptsize}
    \fill (0,0) circle (3.5pt);
\fill (0,2) circle (3.5pt);
\fill (0,-2) circle (3.5pt);

\fill (2,0) circle (3.5pt);
\fill (2,2) circle (3.5pt);
\fill (2,-2) circle (3.5pt);

\fill (-2,0) circle (3.5pt);
\fill (4,0) circle (3.5pt);

\fill (-0.6,3) circle (3.5pt);
\fill (0.6,3) circle (3.5pt);

\fill (1.4,3) circle (3.5pt);
\fill (2.6,3) circle (3.5pt);

\fill (-3,0.6) circle (3.5pt);
\fill (-3,-0.6) circle (3.5pt);

\fill (-0.6,-3) circle (3.5pt);
\fill (0.6,-3) circle (3.5pt);

\fill (1.4,-3) circle (3.5pt);
\fill (2.6,-3) circle (3.5pt);

\fill (5,0.6) circle (3.5pt);
\fill (5,-0.6) circle (3.5pt);
\end{scriptsize}
\draw[color=black] (-0.35,-0.35) node {$v$};
\draw[color=black] (1.6,-0.35) node {$w$};
\end{tikzpicture}
\caption{} \label{FigSharpBothBounds}
\end{figure}
To our best knowledge, Mascia, Rinaldo \cite[Theorem 3.2]{Krull-dim-and-reg-BEI-of-Block-graphs-Rinaldo} give the best lower bound for the regularity of binomial edge ideals of trees. The following theorem gives an improved lower bound for the same.
\begin{theorem}\label{lower bound for reg}
   Let $G$ be a tree that contains $s$ jewels, and $G_{\mathcal J}$ has $p$ connected components. Then  $$\iv(G)+1+D_G-3s-\mu_G+p\leq\reg(S/J_G).$$
\end{theorem}
\begin{proof}
Using Theorem \ref{Gluing theorem}, it is enough to prove the theorem for trees that do not contain a vertex of degree two. In view of Lemma \ref{caterpillar attached to a vertex of a graph}, assume that every caterpillar $T$ with spine $w=w_l,\ldots,w_0$ such that $G=H\cup_w T$ for some $H$, the spine length of $T$ is at most two. We induct on $s$ to prove the theorem. The case $s=1$ follows from Lemma \ref{reg-one-jewel}. Assume $s\geq 2$.  Let $P:v_l,v_{l-1},\ldots,v_0$ be a spine in $G$. By Remark \ref{existance-of-jewel-center-on spine}, $v_2$ is the center of a jewel. Then the graph $G\setminus{\{v_2\}}$ is a collection of $D(v_2)-1$ stars on at least three vertices and $\deg(v_2)-D(v_2)$ isolated vertices. Therefore, by Corollary \ref{reg reln b/w inducedsubgraph}, we have
\begin{equation}\label{lower-bound-delete_v_2 inequality}
    2(D(v_2)-1)+\reg(S/J_{H})\leq \reg(S/J_G),
\end{equation} where $H$ is the non-isolated component of the induced subgraph of $G$ on $(V(G)\setminus{N[v_2]})\cup \{v_3\}$. Using Equation \eqref{lower-bound-delete_v_2 inequality}, it is enough to prove that $\delta-2D(v_2)+2\leq \reg(S/J_H)$, where $\delta=\iv(G)+1+D_G-3s-\mu_G+p$. Note that $H$ is a tree that contains at most $(s-1)$ jewels, and $\iv(H)=\iv(G)-D(v_2)$.

If $D(v_3)\geq 4$, then $H$ contains $(s-1)$ jewels and $D_{H}=D_G-D(v_2)-1$. Note that $H_{\mathcal J}$ has $p$ connected components and $\mu_G=\mu_{H}$. Therefore, by induction, we have $$\delta-2D(v_2)+2\leq \reg(S/J_H).$$ 

Assume $D(v_3)\leq 3$. By Corollary \ref{reg reln b/w inducedsubgraph}, we assume that $\deg_G(v_3)=3$.
Set $$r_1=|\{v\in N_G(v_3)\setminus{\{v_2\}}: \text{$v$ is the center of a jewel, and $D_G(v)=3$}\}|,$$ and $$r_2=|\{v\in N_G(v_3)\setminus{\{v_2\}}: \text{$v$ is the center of a jewel, and $D_G(v)\geq 4$}\}|.$$ Note that $r_1,r_2\in \{0,1,2\}$ and $r_1+r_2\leq 2$. The following table contains all the invariants of $H$ depending on $D(v_3),D(v_4), D(w)$, and $r_1,r_2$.
\newline
\begin{table}[H]
    \centering
    \begin{tabular}{|c|c|c|c|c|c|c|}
    \cline{1-7}
     & & & no. of centers & $D_H$ & $\mu_H$ & $p'$ \\ \cline{1-7}
    \multirow{4}{*}{$D(v_3)=3$} & $r_2=2$ & & \multirow{4}{*}{} & \multirow{4}{*}{} & $\mu_G$ & $p+1$ \\ \cline{2-3}\cline{6-7}
     & \multirow{2}{*}{$r_2=1$} & $D(v_4)\geq 4$ & $s-2-r_1$ & \begin{tabular}{@{}c@{}}
     $D_G-D(v_2)-3$\\
     $-3r_1-r_2$
     \end{tabular}& $\mu_G$ & $p$ \\ \cline{3-3}\cline{6-7}
     & &$D(w)\geq 4$ & &
     & $\mu_H\leq \mu_G+1$ & $\geq p$ \\ \cline{2-3}\cline{6-7}
     & $r_2=0$ & & & & $\mu_H\leq \mu_G+1$ & $\geq p-1$ \\ \cline{1-7}
    \multirow{2}{*}{$D(v_3)\leq 2$} &\multirow{2}{*}{} &$D(v_4)\geq 4$ &$s-1$ &$D_G-D(v_2)-1$ &$\mu_G-1$ &$p-1$ \\ \cline{3-7}
     & & $D(v_4)=3$ &$s-2$ &$D_G-D(v_2)-3$ &$\leq \mu_G$ &$\geq p-1$ \\ \cline{3-7}
      & & $D(v_4)\leq 2$ &$s-1$ &$D_G-D(v_2)$ &$\mu_G$ &$ p-1$ \\ \cline{1-7}
    \end{tabular}
    \caption{}
    \label{tab:placeholder}
\end{table}
$\textbf{Case 1}:$ Suppose $D(v_3)=3$.
Then $H$ contains $s-2-r_1$ centers with $D_H=D_G-D(v_2)-3-3r_1-r_2$, and  ${H}_{\mathcal
J}$ has at least $p-1$ connected components. Note that $C_H\setminus{C_G}\subset\{v_4\}$. Therefore, $\mu_H\leq \mu_G+1$. More precisely, if $v_2,v_4$ and $w$ are supporting vertices of $\mathcal J_{v_3}$, and $p'$ is the number of the connected components of $H_{\mathcal{J}}$, then 
\begin{enumerate}
    \item if $r_2=2$, then $p'=p+1$ and $\mu_H=\mu_G$,
    \item if $r_2=1$, then 
    \begin{enumerate}
        \item if $D(v_4)\geq 4$, then $p'= p$ and $\mu_H=\mu_G$,
        \item if $D(w)\geq 4$, then $p'\geq p$. In fact, if $\mu_H=\mu_G+1$, then $p'= p+2$.
    \end{enumerate}
    \item if $r_2=0$, then $p'\geq p-1$. In fact, if $\mu_H=\mu_G+1$, then $p'=p+1$. 
    
\end{enumerate}
Thus, using induction, we have $\delta-2D(v_2)+2\leq \reg(S/J_{H})$.
\newline
$\textbf{Case 2}:$ Suppose $D(v_3)\leq 2$. If $D(v_4)\geq 4$, then $H$ contains $(s-1)$ jewels with $D_H=D_G-D(v_2)-1$. The number of connected components of $H_{\mathcal J}$ is $p-1$, and $\mu_H= \mu_G-1$. If $D(v_4)\leq2 $, then $H$ contains $(s-1)$ jewels with $D_H=D_G-D(v_2)$. In this case, $p'=p-1$ and $\mu_H=\mu_G$. If $D(v_4)=3$, then $H$ contains $(s-2)$ jewels with $D_H=D_G-D(v_2)-3$. Further, we have $\mu_H\leq \mu_G$ and  $p'\geq p-1$. In either case, we obtain the desired result by induction. 
\end{proof} 
\begin{example}
    Let $G$ be the graph in Figure \ref{FigSharpBothBounds}. As noted in Example \ref{ex: upper bound}, we have $\iv(G)=8,s=2,D_G=8$, and $p=1$. Moreover, there is no vertex $u\in V(G)$ with $D(u)=2$, and hence $\mu_G=0$. Therefore, by Theorem \ref{lower bound for reg}, we obtain $\reg(S/J_G)\geq 12$.
\end{example}
We arrive at Corollary \ref{Trees of exact reg} that gives the regularity of binomial edge ideals of trees under some condition. 
\begin{corollary}\label{Trees of exact reg}
   Let $G$ be a tree containing $s$ jewels with $\mu_G=0$ in Theorem \ref{lower bound for reg} and $s_i\leq 2$ for $1\leq i\leq p$. Then $$\reg(S/J_G)=\iv(G)+1+D_G-2s-e_G.$$ 
\end{corollary}
\begin{proof}
    Since $\mu_G=0$, the lower bound of Theorem \ref{lower bound for reg} is $\iv(G)+1+D_G-3s+p$. Using the fact that $s_i\leq 2$ for all $1\leq i\leq p$ and $\sum\limits_{i=1}^ps_i=s$, we get $3s-p=2s+\sum\limits_{i=1}^{p}(s_i-1)=2s+e_G$. Hence, $\reg(S/J_G)\geq\iv(G)+1+D_G-2s-e_G$, and the result follows from Theorem \ref{upper bound for reg}.
\end{proof}
\begin{remark}\label{Remark-other bounds are strict}
\label{ExRegEquality}
   Let $G$ be the tree shown in Figure \ref{FigSharpBothBounds}. Then $\reg(S/J_G)=12$, which coincides with the lower bound in Theorem \ref{lower bound for reg} and the upper bound in Theorem \ref{upper bound for reg}. However, the lower bound of \cite[Corollary 3.2]{Krull-dim-and-reg-BEI-of-Block-graphs-Rinaldo}, and the upper bound of \cite[Corollary 4.7]{BEI-of--generalized-block-graphs-Arvind-Kumar} are strict. 
\end{remark} 
\begin{example}\label{Exm bound better than kumar}
 Let $G$ be the graph obtained from $d$ copies of $\Gamma_v$ shown in Figure \ref{FigBoundIsBetterThanKumar} by identifying the vertex $v$, where $d\geq 3$. Then $\iv(G)=3d+1$. Since the adjacent vertex to $v$ of $\Gamma_v$, and $v$ itself, becomes the center of a jewel in $G$, it follows that $s=d+1$. Since $D_G(v)=d$ and every other center $w$ satisfies $D_G(w)=3$, it follows that $D_G=4d$. Moreover, the graph $G_{\mathcal J}$ is $K_{1,d}$, and hence $p=1$ with $s_1=p+1$ and $e_G=0$. By Theorem \ref{upper bound for reg}, we obtain $\reg(S/J_G)\leq 5d-\left\lfloor \frac{d+1}{3}\right\rfloor$, while the bound in \cite[Corollary 4.7]{BEI-of--generalized-block-graphs-Arvind-Kumar} is $5d$. Notice that as $d$ increases, the difference between the bounds of Theorem \ref{upper bound for reg} and that of \cite[Corollary 4.7]{BEI-of--generalized-block-graphs-Arvind-Kumar} also increases.
\end{example}
\begin{figure}[h]
\centering
    \begin{tikzpicture}[scale=0.4, line cap=round,line join=round,x=1cm,y=1cm]
\clip(-11,-6) rectangle (10,6);
\draw (-7,0)-- (-7,2);
\draw (-7,2)-- (-7.6,3);
\draw (-7,2)-- (-6.4,3);
\draw (-7,0)-- (-9,0);
\draw (-9,0)-- (-10,0.6);
\draw (-9,0)-- (-10,-0.6);
\draw (-7,0)-- (-7,-2);
\draw[color=black] (-7.5,-2.5) node {$v$};
\draw[color=black] (-7,-5) node {$\Gamma_{v}$};
\draw  (0,0)-- (0,2);
\draw  (0,2)-- (-0.6,3);
\draw  (0,2)-- (0.6,3);
\draw  (0,0)-- (-2,0);
\draw  (-2,0)-- (-3,0.6);
\draw (-2,0)-- (-3,-0.6);
\draw  (0,0)-- (2,0);
\draw  (2,0)-- (2.84,1.76);
\draw  (2.84,1.76)-- (2.82,3.56);
\draw  (2.84,1.76)-- (4.82,1.8);
\draw  (2.82,3.56)-- (2.24,4.56);
\draw  (2.82,3.56)-- (3.42,4.58);
\draw  (4.82,1.8)-- (5.8,2.38);
\draw  (4.82,1.8)-- (5.8,1.18);
\draw  (2,0)-- (3.58,-1.38);
\draw  (3.58,-1.38)-- (5.6,-1.44);
\draw  (3.58,-1.38)-- (3.62,-3.42);
\draw  (5.6,-1.44)-- (6.6,-0.82);
\draw  (5.6,-1.44)-- (6.6,-2);
\draw  (3.62,-3.42)-- (4.22,-4.4);
\draw  (3.62,-3.42)-- (3,-4.42);
\draw (-1.12,-0.67);
\draw (-0.74,-1.94) circle (1pt);
\draw (-0.6,-1.57);
\draw (0.06,-2.82) circle (1pt);
\draw (0.25,-2.4) ;
\draw (.9,-3.40) circle (1pt);
\draw (1.19,-2.82) ;
\draw[color=black] (1,-5) node {$G$};
\draw[color=black] (1.65,-.4) node {$v$};
\begin{scriptsize}
    \fill (-7,0) circle (3.5pt);
    \fill (-7,2) circle (3.5pt);
\fill (-7,-2) circle (3.5pt);
\fill (-7.6,3) circle (3.5pt);
\fill (-6.4,3) circle (3.5pt);

\fill (-9,0) circle (3.5pt);
\fill (-10,0.6) circle (3.5pt);
\fill (-10,-0.6) circle (3.5pt);

\fill (0,0) circle (3.5pt);
\fill (0,2) circle (3.5pt);
\fill (-0.6,3) circle (3.5pt);
\fill (0.6,3) circle (3.5pt);

\fill (-2,0) circle (3.5pt);
\fill (-3,0.6) circle (3.5pt);
\fill (-3,-0.6) circle (3.5pt);

\fill (2,0) circle (3.5pt);
\fill (2.84,1.76) circle (3.5pt);
\fill (2.82,3.56) circle (3.5pt);
\fill (2.24,4.56) circle (3.5pt);
\fill (3.42,4.58) circle (3.5pt);

\fill (4.82,1.8) circle (3.5pt);
\fill (5.8,2.38) circle (3.5pt);
\fill (5.8,1.18) circle (3.5pt);

\fill (3.58,-1.38) circle (3.5pt);
\fill (5.6,-1.44) circle (3.5pt);
\fill (6.6,-0.82) circle (3.5pt);
\fill (6.6,-2) circle (3.5pt);

\fill (3.62,-3.42) circle (3.5pt);
\fill (4.22,-4.4) circle (3.5pt);
\fill (3,-4.42) circle (3.5pt);
\end{scriptsize}
\end{tikzpicture}
\caption{} \label{FigBoundIsBetterThanKumar}
\end{figure}
 \subsection*{Data availability statement} Data sharing is not applicable to this article as no new data were created or analyzed in this study.

 \subsection*{Conflict of interest} The authors declare that they have no known competing financial interests or personal relationships that could have appeared to influence the work reported in this paper.
\bibliographystyle{plain}
	\bibliography{bib}
\end{document}